\documentclass[12pt]{article}
\usepackage{amsmath, amssymb}

\newcommand{\qed}{$\hfill \Box$}
\newcommand{\proof}{\medskip \noindent {\bf Proof:}\qquad}

\newtheorem{theorem}{Theorem}[section]

\newtheorem{proposition}{Proposition}[section]

\newcommand{\EE}{{\mathbb E}}
\newcommand{\KP}{{\cal K}^{\cal P}}
\newcommand{\cS}{{S}}

\newlength{\IndentI}
\newlength{\IndentII}
\newlength{\IndentIII}
\newlength{\WidthI}
\newlength{\WidthII}
\newlength{\WidthIII}
\title{Game theoretic derivation of discrete distributions and discrete
  pricing formulas}
\author{
Akimichi Takemura and Taiji Suzuki\\
Department of Mathematical Informatics\\
Graduate School of Information Science and Technology\\
University of Tokyo\\
}

\date{September, 2005}

\begin{document}
\maketitle
\begin{abstract}
  In this expository paper we illustrate the generality of game theoretic
  probability protocols of Shafer and Vovk (2001) in finite-horizon
  discrete games. By restricting ourselves to finite-horizon discrete
  games, we can explicitly describe how discrete distributions with finite
  support and the discrete pricing formulas, such as the
  Cox-Ross-Rubinstein formula, are naturally derived from game-theoretic
  probability protocols.  Corresponding to any discrete distribution with
  finite support, we construct a finite-horizon discrete game, a
  replicating strategy of Skeptic, and a neutral forecasting strategy of
  Forecaster, such that the discrete distribution is derived from the
  game.  Construction of a replicating strategy is the same as in the
  standard arbitrage arguments of pricing European options in the binomial
  tree models.  However the game theoretic framework is advantageous
  because no a priori probabilistic assumption is needed.
\end{abstract}


\section{Introduction}
\label{sec:intro}

In the game theoretic probability of Shafer and Vovk (2001), probability
distributions and probability models are not assumed a priori but derived
as logical consequences of certain protocol of a game between two players
``Skeptic'' and ``Reality''.  In this game Skeptic tries to become rich by
exploiting patterns in the moves of Reality.  In order to prevent Skeptic
from becoming rich, Reality is in a sense forced to behave
probabilistically.  Therefore probability distributions are determined by
the protocol of the game.  This feature of the game theoretic probability
is well illustrated by Shafer and Vovk (2001) in their derivation of
Skeptic's strategy forcing the strong law of large numbers (Chapter 3)
and the derivation of Black-Scholes formula (Chapter 9).
Also in Takeuchi's exposition of the game theoretic probability and
finance (Takeuchi (2004)) this point is discussed with many interesting
examples.  Recently Kumon and Takemura (2005) gave a very simple strategy
forcing the strong law of large numbers.

In the standard stochastic derivation of option pricing formulas,
empirical probability is assumed first, but then by arbitrage arguments,
the empirical probability is replaced by the risk neutral probability and
the price of an option is given as the expected value with respect to the
risk neutral probability.  The risk neutral probability is often
explained as a purely operational device useful in expressing the option
price in a convenient form.  On the other hand in the game theoretic
probability the risk neutral probability is more substantial, in the sense
that Reality is forced to behave according to the risk neutral probability
to avoid arbitrage by Skeptic.  We should mention here that in Shafer and Vovk
(2001) ``forcing'' is used only for infinite-horizon
games.  
In this paper we somewhat informally use the word to mean that Reality
should avoid arbitrage by Skeptic in the setting of finite-horizon games.

Additional flexibility of game theoretic probability is gained by
introducing the third player ``Forecaster'' into the game.  At the
beginning of each round Forecaster sets the price for Reality's move.  By
appropriately specifying the strategy of Forecaster, Reality's moves can
be forced to follow any prespecified distribution.

In this paper we demonstrate the above features of the game theoretic
probability in the setting of finite-horizon discrete games.  For
expository purposes we start with the simplest setting of the coin-tossing
game and derive binomial distribution in Section \ref{sec:binomial} and
give an  analogous derivation of the Cox-Ross-Rubinstein formula in Section
\ref{sec:cox-ross-rubinstein}.  We discuss derivation of hypergeometric
distribution and Polya's distribution in Section \ref{sec:hypergeometric}
in order to illustrate the role of Forecaster.  Then in Section
\ref{sec:arbitrary-distribution} we discuss derivation of an arbitrary discrete
distribution with finite support.  Multivariate extension is given in
Section \ref{sec:multivariate}.
Some preliminary material on game theoretic
probability is given in Appendix.

%
%

\section{Derivation of binomial distribution}
\label{sec:binomial}

Consider the finite-horizon fair-coin game in  
Section 6.1 of Shafer and Vovk (2001).
Its protocol is given as
follows.

\medskip
\noindent
\textsc{Fair-Coin Game}\\
\textbf{Protocol:}

\parshape=6
\IndentI   \WidthI
\IndentI   \WidthI
\IndentII  \WidthII
\IndentII  \WidthII
\IndentII  \WidthII
\IndentI   \WidthI
\noindent
${\cal K}_0 =\alpha$ : given \\
FOR $n=1,\dots,N$\\
  Skeptic announces $M_n\in{\mathbb R}$.\\
  Reality announces $x_n\in \{-1,1\}$.\\
  ${\cal K}_n := {\cal K}_{n-1} + M_n x_n$.\\
END FOR

\medskip

In this protocol a game theoretic version of Chebyshev  inequality is 
proved in (6.9) of S\&V  in the following form:
\[
\underline{\mathbb P} \left\{ \left| \frac{\cS_N}{N} \right |  \le
    \epsilon
\right\}  \ge 1 - \frac{1}{N \epsilon^2},
\]
where $\cS_N = x_1 + \dots + x_N$ and $\underline{\mathbb P}$ denotes the
lower probability.  Actually the
equality of the upper probability and the lower probability 
\begin{equation}
\label{eq:chebyshev}
\underline{\mathbb P} \left\{ \left| \frac{\cS_N}{N} \right |  \le
    \epsilon
    \right\}  
=
 \bar{\mathbb P} \left\{ \left| \frac{\cS_N}{N} \right |  \le
    \epsilon
    \right\} 
\end{equation}
holds here and this probability is given by binomial distribution.
Although this fact is contained in a more general statement of Proposition
8.5 of S\&V, we give a full proof of this fact employing standard arbitrage
arguments.

In order to treat success probability $p\neq 1/2$, let us consider
the following biased-coin game.

\medskip
\noindent
\textsc{Biased-Coin Game}\\
\textbf{Protocol:}

\parshape=6
\IndentI   \WidthI
\IndentI   \WidthI
\IndentII  \WidthII
\IndentII  \WidthII
\IndentII  \WidthII
\IndentI   \WidthI
\noindent
${\cal K}_0 =\alpha\in {\mathbb R}, \ a,b>0$: given\\
FOR $n=1,\dots,N$\\
  Skeptic announces $M_n\in{\mathbb R}$.\\
  Reality announces $x_n\in \{a, -b\}$.\\
  ${\cal K}_n := {\cal K}_{n-1} + M_n x_n$.\\
END FOR\\
\noindent

As above we write $\cS_n = x_1 + \dots + x_n$. Let 
$\cS_0=0$. Consider a random variable $x(\xi)=\eta(\cS_N)$ which
depends only on $\cS_N$ (European option).   Then we have the
following basic result.

\begin{theorem}
\label{thm:binomial}
  The upper and the lower expected values of $\eta(\cS_N)$ coincide
  and given by
\begin{equation}
\label{eq:eesn}
\bar\EE(\eta(\cS_N)) = \underline\EE(\eta(\cS_N)) = 
\sum_{m=0}^N {N \choose m} p^m (1-p)^{N-m} \eta(ma - (N-m)b),
\end{equation}
where $p=b/(a+b)$ is the risk neutral probability.
\end{theorem}

\proof The first step of our proof consists of defining a
``candidate'' price of the European option.
In the second step we verify that the candidate price is actually the precise
price, by constructing a replicating strategy.

Let $\bar \eta(n,\cS_n)$, $\cS_n= -nb, -(n-1)b+a, \ldots na$,
denote the price of $\eta(\cS_N)$ at time $n$.  We require 
$\eta(n,\cS_n)$  to 
satisfy  the following ``partial difference equation''
\begin{equation}
\label{eq:heat-discrete}
\bar \eta(n, \cS_n)=  p \bar \eta(n+1,\cS_n + a) + q \bar \eta(n+1,
  \cS_n -b), \quad  0 \le n < N, 
\end{equation}
where $q=1-p$.  Note that (\ref{eq:heat-discrete}) with $p=1/2$
is a discrete version of the heat equation. 
The terminal condition for $\bar \eta(n, \cS_n)$ is given by
\begin{equation}
\label{eq:discrete-terminal}
\bar \eta(N, \cS_N) = \eta(\cS_N),\quad \cS_N=ma - (N-m)b, \ m=0,\dots,N.
\end{equation}
Starting with the terminal condition (\ref{eq:discrete-terminal}) we can
solve for $\bar \eta(n, \cS_n)$ in 
(\ref{eq:heat-discrete}) by backward induction $n=N-1, N-2, \dots,0$.
Then the initial value
$\bar \eta(0,0)$ is easily calculated as
\begin{eqnarray*}
\bar \eta(0,0) &=& p \bar \eta(1,a) + q \bar \eta(1,-b) \\
  &=& p( p\bar \eta(2,2a) + q \bar \eta(2,a-b))  
     + q( p\bar \eta(2,a-b)+ q \bar \eta(2, -2b))\\
  &=& p^2 \bar \eta(2,2a) + 2pq \bar \eta(2,a-b) + q^2 \bar \eta(2,-2b)\\
  &=& \dots\\
  &=& \sum_{m=0}^N {N \choose m} p^m (1-p)^{N-m} \eta(ma - (N-m)b).
\end{eqnarray*}

Now we  describe a replicating strategy for
$\eta(\cS_N)$ with the the replicating initial capital $\bar \eta(0,0)$.  
For $n=1,\dots,N$, let
\begin{equation}
\label{eq:mn1}
M_n = \frac{\bar \eta(n,\cS_{n-1}+a) - \bar \eta(n, \cS_{n-1}-b)}{a+b} .
\end{equation}
Note that $a+b$ can be written as 
$$
a+b=(\cS_{n-1}+a) - (\cS_{n-1}-b).
$$
Therefore $M_n$ is the ratio of the increments of $\bar\eta(n,{\cal
  S}_n)$ and $\cS_n$ and is called the ``delta hedge''.
We now check that this $M_n$ gives a replicating strategy $\cal P$.  This can be
confirmed by forward induction.  At the end of the first round $n=1$, 
\begin{eqnarray*}
  \bar\eta(0,0) + {\cal K}^{\cal P}_1 
  &=& \bar \eta(0,0) + M_1 x_1\\
  &=& \bar \eta(0,0) + \frac{\bar \eta(1,a)-\bar \eta(1,-b)}{a+b} x_1\\
  &=& 
  \begin{cases} 
    \bar \eta(0,0) + q (\bar \eta(1,a) - \bar \eta(1, -b)), & {\rm if}\   x_1 = a, \\
    \bar \eta(0,0) - p (\bar \eta(1,a) - \bar \eta(1, -b)), & {\rm if}\
    x_1 = -b
  \end{cases}\\
  &= &
\begin{cases} 
  \bar \eta(1,a), & {\rm if}\   x_1 = a, \\
  \bar \eta(1,-b) , & {\rm if}\ x_1 = -b,
\end{cases}\\
&=& \bar \eta(1,\cS_1).
\end{eqnarray*}
Similarly at the end of round $n=2$, we have
\begin{eqnarray*} \bar\eta(0,0)+ {\cal K}^{\cal P}_2 
&=&
\begin{cases} 
  \bar \eta(2,\cS_1 + a), & {\rm if}\  x_2 = a, \\
  \bar \eta(2,\cS_1 - b), & {\rm if}\  x_2 = -b,
\end{cases}  \\
&=&
\bar \eta(2,\cS_2).
\end{eqnarray*}
Now by induction we arrive at
\[
\bar\eta(0,0)+{\cal K}^{\cal P}_N = \bar \eta(N,\cS_N)=\eta(\cS_N).
\]
We have confirmed that $M_n$  in (\ref{eq:mn1}) with the replicating
initial capital (\ref{eq:eesn}) gives a replicating strategy for
$\eta(\cS_N)$.  Hence the theorem holds by Proposition
\ref{prop:replicatioting-strategy} in  Appendix.  \qed

\bigskip
In particular if we take
\[
\eta(\cS_N)= 
\begin{cases} 
  1, & {\rm if}\   |\cS_N|/N \le \epsilon , \\
  0, & {\rm otherwise},
\end{cases} 
\]
we see that the equality holds in  (\ref{eq:chebyshev}) and the
probability is given by binomial distribution.

In this section we took Reality's move space as $\{a,-b\}$.  This is
convenient in comparing Theorem \ref{thm:binomial} with the
Cox-Ross-Rubinstein formula in the next section.  However for generalization
of binomial distribution to hypergeometric distribution in Section
\ref{sec:hypergeometric},
it is more convenient to rescale Reality's move space
to $\{ 0,1\}$.  Then we need to introduce the price $p$ for the ``ticket''
$x_n$.  The rescaled protocol is written as follows.

\medskip
\noindent
\textsc{Rescaled Biased-Coin Game}\\
\textbf{Protocol:}

\parshape=6
\IndentI   \WidthI
\IndentI   \WidthI
\IndentII  \WidthII
\IndentII  \WidthII
\IndentII  \WidthII
\IndentI   \WidthI
\noindent
${\cal K}_0 =\alpha, 0<p<1$: given\\
FOR $n=1,\dots,N$\\
  Skeptic announces $M_n\in{\mathbb R}$.\\
  Reality announces $x_n\in \{0, 1\}$.\\
  ${\cal K}_n := {\cal K}_{n-1} + M_n (x_n-p)$.\\
END FOR\\
\noindent

It is clear that the biased-coin game and the rescaled biased-coin game is
equivalent by the affine correspondence
$x_n \leftrightarrow (a+b)(x_n-p),  p=b/(a+b)
$.  In the rescaled version 
the expected value in (\ref{eq:eesn}) is simply written as
\begin{equation}
\label{eq:eesnr}
\bar\EE(\eta(\cS_N)) = \underline\EE(\eta(\cS_N)) = 
\sum_{m=0}^N {N \choose m} p^m (1-p)^{N-m} \eta(m).
\end{equation}
Furthermore, since the increment of $\cS_n$ is normalized to be 1, 
the replicating strategy in (\ref{eq:mn1}) is simply written as
\begin{equation}
\label{eq:mn1s}
M_n = \bar \eta(n,\cS_{n-1}+1) - \bar \eta(n, {\cS}_{n-1}).
\end{equation}

It is also conceptually very important to consider the single step game
i.e.\ the game with $N=1$. Note that each round $n$ of the $N$ step
biased-coin game can be viewed as a single step game. In the single step
game binomial distribution reduces to a Bernoulli trial.  This implies
that given the price $p$, Reality's move $x_n$ for each round $n$ is
exactly the same as a single Bernoulli trial with success probability $p$.
Furthermore this behavior of Reality is dictated solely by the value of
$p$, independently from the past moves $x_1,\ldots,x_{n-1}$ of Reality.
Therefore in the Rescaled Biased-Coin Game, Reality's moves $x_1,
\dots,x_N$ are independent Bernoulli trials.

\section{Derivation of the Cox-Ross-Rubinstein formula}
\label{sec:cox-ross-rubinstein}

Here we present a game theoretic formulation and derivation of the
Cox-Ross-Rubinstein formula (Cox, Ross and Rubinstein (1979)), which is
fully discussed in many introductory textbooks on option pricing (e.g.,
Shreve (2004), Chapter 2 of Baxter and Rennie (1996), Chapter 8 of Capi\'nski and
Zastawniak (2003)).  Once an appropriate game is formulated, the rest of
the argument is the same as in the previous section.

Our protocol for Cox-Ross-Rubinstein game is as follows.

\medskip
\noindent
\textsc{Cox-Ross-Rubinstein Game}\\
\textbf{Protocol:}

\parshape=7
\IndentI   \WidthI
\IndentI   \WidthI
\IndentII  \WidthII
\IndentII  \WidthII
\IndentII  \WidthII
\IndentII  \WidthII
\IndentI   \WidthI
\noindent
$\cS_0 > 0,   u>r>d > 0$: given\\
FOR $n=1,\dots,N$\\
  Skeptic announces $M_n\in{\mathbb R}$.\\
  Reality announces $x_n\in \{u, d\}$.\\
  $\cS_n:=\cS_{n-1} \times x_n$.\\
  ${\cal K}_n := {\cal K}_{n-1} + M_n (\cS_n - \cS_{n-1})
   + (r-1) ({\cal K}_{n-1}- M_n \cS_{n-1})$.\\
END FOR\\
\noindent

Here ${\cal K}_{n-1}- M_n \cS_{n-1}$ is the amount of riskless bond
held by Skeptic for the round $n$ and $r-1$ is the fixed riskless interest
rate.  Although by appropriate discounting we may put $r=1$ without
essential loss of generality (Section 12.1 of S\&V), here we leave the
interest rate $r$ as in standard derivation of the
Cox-Ross-Rubinstein formula.  $p= (r-d)/(u-d)$ is called the risk
neutral probability.

Let $\eta(\cS_N)$ denote a payoff function of a European option
depending on $\cS_N$.  Corresponding to Theorem \ref{thm:binomial} we
have the following result.

\begin{theorem} (the Cox-Ross-Rubinstein formula)\qquad 
\label{thm:crr}
  The upper and the lower expected values of $\eta(\cS_N)$ coincide
  and given by
\begin{equation}
\label{eq:crrf}
\bar\EE(\eta(\cS_N)) = \underline\EE(\eta(\cS_N)) = 
\frac{1}{r^N} \sum_{m=0}^N {N \choose m} p^m (1-p)^{N-m} \eta(u^m d^{N-m}
\cS_0),
\end{equation}
where $p=(r-d)/(u-d)$ is the risk neutral probability.
\end{theorem}

\proof As in the previous section we define $\bar \eta(n,\cS_n)$ by
backward induction. Put 
$\bar \eta(N,\cS_N)=\eta(S_N)$ and for  $n=N-1,\dots,0$, define
\[
\bar \eta(n,\cS_n) = \frac{1}{r} \Big(p \bar \eta(n+1, u\cS_n) + (1-p)
  \bar \eta(n+1, d \cS_n)\Big)
\]
Then the initial value $\bar \eta(0,\cS_0)$ is easily calculated as
\[
\bar \eta(0,\cS_0) = \frac{1}{r^N} \sum_{m=0}^N {N \choose m} p^m
(1-p)^{N-m}  \eta(u^m d^{N-m}\cS_0).
\]
This becomes the  replicating initial capital of the following replicating
strategy:
\[
M_n = \frac{\bar\eta(n,u \cS_{n-1}) - \bar \eta(n, d\cS_{n-1})}
      {(u-d) \cS_{n-1}}.
\]
Since the game is coherent by the requirement $u > r > d$, 
the theorem follows from Proposition
\ref{prop:replicatioting-strategy}.
\qed

\section{Hypergeometric distribution and Polya's distribution}
\label{sec:hypergeometric}

In the rescaled biased-coin game of Section \ref{sec:binomial}, the price
$p$ of the ticket $x_n$ was a constant.  Therefore the third player
``Forecaster'' did not enter the protocol.  Now we introduce Forecaster,
who sets the price of the ticket at the beginning of each round in the
rescaled biased-coin game.  We illustrate the role of Forecaster below by
deriving the hypergeometric distribution.  We also derive Polya's distribution.

\medskip
\noindent
\textsc{Rescaled Biased-Coin Game With Forecaster}\\
\textbf{Protocol:}

\parshape=7
\IndentI   \WidthI
\IndentI   \WidthI
\IndentII  \WidthII
\IndentII  \WidthII
\IndentII  \WidthII
\IndentII  \WidthII
\IndentI   \WidthI
\noindent
${\cal K}_0 =\alpha$: given\\
FOR $n=1,\dots,N$\\
  Forecaster announces $p_n \in {\mathbb R}$.\\
  Skeptic announces $M_n\in{\mathbb R}$.\\
  Reality announces $x_n\in \{0, 1\}$.\\
  ${\cal K}_n := {\cal K}_{n-1} + M_n (x_n-p_n)$.\\
END FOR\\
\noindent

Note that if Forecaster announces $p_n > 1$ or
$p_n < 0$, then Skeptic can become infinitely rich immediately by taking $|M_n|$
arbitrarily large.  Therefore we can restrict the move space of Forecaster
to be $[0,1]$.  Furthermore if $p_n=0$,  Skeptic can still take $M_n$
arbitrarily large, which forces Reality to choose $x_n=0$.  Similarly if
$p_n=1$, then Reality is forced to choose $x_n=1$.

Now consider a strategy of Forecaster.  A strategy of Forecaster is called
{\it neutral\/} (Section 8.2 of S\&V) if $p_n$ is determined by the past
moves of Reality $x_1\ldots x_{n-1}$.  {}From now on we only consider neutral
strategies for Forecaster.  Furthermore for simplicity we consider
neutral strategy depending on $\cS_{n-1}=x_1+\cdots+x_{n-1}$ and write
$p_n=p_n(\cS_{n-1})$, which we may call ``Markovian neutral strategy''.
In Markovian neutral strategy Forecaster only needs to keep $\cS_{n-1}$ in
memory to choose his move.

Consider an urn with 
$\nu_1$ red balls and $\nu_2$ black
balls, where $\nu_1 + \nu_2 \ge N$. 
Let $x_n=1$  correspond to drawing a red ball and let
$x_n=0$ correspond to drawing a black ball from the urn by Reality.
Let $p_n$ be the ratio of red balls in the urn at the $n$-th round. Then
$$
p_n = p_n(\cS_{n-1})=
  \frac{\max(0,\nu_1 - \cS_{n-1})}{\nu_1 + \nu_2 - (n-1)}. 
$$
Actually here we do need to take the positive part of $\nu_1 - \cS_{n-1}$,
because as remarked above once the boundary $\cS_{n-1}=\nu_1$ is attained,
then $0=p_n=p_{n+1}=\dots$ and Reality is forced to choose
$0=x_n=x_{n+1}=\dots$, which results in $\nu_1 = \cS_n = \cS_{n+1}=\dots$.
Now we write out a game of sampling without replacement from an urn.

\medskip
\noindent
\textsc{Game of Sampling Without Replacement From An Urn}\\
\textbf{Protocol:}

\parshape=8
\IndentI   \WidthI
\IndentI   \WidthI
\IndentII  \WidthII
\IndentII  \WidthII
\IndentII  \WidthII
\IndentII  \WidthII
\IndentII  \WidthII
\IndentI   \WidthI
\noindent
${\cal K}_0 =\alpha, \nu_1>0, \nu_2>0, \nu_1+\nu_2 \ge N$, 
$\cS_0=0$  : given\\
FOR $n=1,\dots,N$\\
  Forecaster announces $p_n =(\nu_1-\cS_{n-1})/(\nu_1 + \nu_2 - n+1)$\\
  Skeptic announces $M_n\in{\mathbb R}$.\\
  Reality announces $x_n\in \{0, 1\}$.\\
  ${\cal K}_n := {\cal K}_{n-1} + M_n (x_n-p_n)$.\\
  $\cS_n := \cS_{n-1}+x_n$.\\
END FOR\\
\noindent

For this game the upper and the lower values of the payoff $\eta({\cal
  S}_N)$  coincide and are given  by the expected value with respect to
the hypergeometric distribution:
\begin{equation}
\label{eq:eehy}
\bar\EE(\eta(\cS_N)) = \underline\EE(\eta(\cS_N)) = 
\sum_{m=\max(0,N-\nu_2)}^{\min(\nu_1,N)} \eta(m) 
\frac{ {\nu_1 \choose m} {\nu_2 \choose {N-m} } } 
{ {{\nu_1 + \nu_2}  \choose  N} } .
\end{equation}
This result is actually almost obvious from the discussion at the end of
Section \ref{sec:binomial}, namely, at each round $n$ Reality's move $x_n$
is like drawing a ball from an urn with $\nu_1 - \cS_{n-1}$ red balls and
$\nu_2 - (n-1- \cS_{n-1})$ black balls.  However it is instructive to look
at a formal proof of (\ref{eq:eehy}).

Define a candidate price of $\eta(\cS_N)$ at time $n$ by
backward induction:
\[
  \bar \eta(n,\cS_n)= p_{n+1}(\cS_n)   \bar\eta(n+1,\cS_n+1)
  + (1-p_{n+1}(\cS_n)) \bar\eta(n+1,\cS_n) , \ n=N-1, \ldots,0,
\]
with the terminal condition $\bar \eta(N,\cS_N)= \eta(\cS_N)$. Then by fully
expanding the recurrence relation we have
\begin{equation}
\label{eq:full-expansion}
\bar\eta(0,0)=\sum_{(x_1, \ldots, x_N)\in { \{0,1\}^N}\atop
 \max(0,N-\nu_2)\le \cS_N \le \min(\nu_1,N)}
\prod_{n=1}^N   p_n(\cS_{n-1})^{x_n} (1-p_n(\cS_{n-1}))^{1-x_n}  
\eta(\cS_N).
\end{equation}
Actually we do not need the restriction 
$ \max(0,N-\nu_2)\le \cS_N \le \min(\nu_1,N)$ in the summation, because
$\prod_{n=1}^N   p_n(\cS_{n-1})^{x_n} (1-p_n(\cS_{n-1}))^{1-x_n} =0$ 
for $\cS_N$ outside of this range.
Now it is easily seen that
\begin{eqnarray*}
&& \prod_{n=1}^N   p_n(\cS_{n-1})^{x_n} (1-p_n(\cS_{n-1}))^{1-x_n} \\
&& \qquad \qquad = 
  \frac{\nu_1 (\nu_1-1)\cdots (\nu_1-\cS_N+1)\; \cdot \;
  \nu_2 (\nu_2-1)\cdots (\nu_2-N+\cS_N+1)}
{(\nu_1 + \nu_2) (\nu_1 + \nu_2 - 1) \cdots (\nu_1 + \nu_2 - N+1)}\\
&& \qquad\qquad 
= \frac{\frac{\nu_1!}{(\nu_1-{\cS}_N)!} \frac{\nu_2!}{(\nu_2-(N-{\cS}_N))!}}
{\frac {(\nu_1 + \nu_2)!} {(\nu_1 + \nu_2 -N)!}} .
\end{eqnarray*}
Therefore
\[
\bar\eta(0,0)=\sum_{(x_1, \ldots, x_N)\in { \{0,1\}^N}
 \atop  \max(0,N-\nu_2)\le \cS_N \le \min(\nu_1,N)}
\frac{\frac{\nu_1!}{(\nu_1-{\cS}_N)!} \frac{\nu_2!}{(\nu_2-(N-{\cS}_N))!}}
{\frac {(\nu_1 + \nu_2)!} {(\nu_1 + \nu_2 -N)!}}
\eta(\cS_N).
\]
Now for a given value of $\cS_N$, the summation just counts the number of
ways of choosing $\cS_N$ 1's among $x_1, \ldots, x_N$.  It follows that
\[
\bar\eta(0,0)=\sum_{\max(0,N-\nu_2)\le \cS_N \le \min(\nu_1,N)}
\frac{N!}{\cS_N! (N-\cS_N)!} 
\frac{\frac{\nu_1!}{(\nu_1-{\cS}_N)!} \frac{\nu_2!}{(\nu_2-(N-{\cS}_N))!}}
{\frac {(\nu_1 + \nu_2)!} {(\nu_1 + \nu_2 -N)!}}
\eta(\cS_N),
\]
which proves (\ref{eq:eehy}).

The above argument can be immediately applied to Polya's urn model
(Section V.2 of Feller (1968)).  In this scheme, when a ball is drawn from
an urn, it is replaced and, moreover, $c$ balls of the same color 
are added.  Then the game corresponding to
Polya's urn model differs from the game of sampling without replacement
only in the specification of Forecaster's neutral strategy.  In Polya's
urn model
\[
p_n = p_n(\cS_{n-1}) 
= \frac{\nu_1 + c \cS_{n-1}}{ \nu_1 + \nu_2 + (n-1)c}.
\]
Then as in (2.4) of Section V.2 of Feller (1968), the expected value of
$\eta({\cS}_N)$ in this game is written as
\begin{equation}
\label{eq:eepolya}
\bar\EE(\eta(\cS_N)) = \underline\EE(\eta(\cS_N)) = 
\sum_{m=0}^N \eta(m) 
\frac{ {-\nu_1/c \choose m} {-\nu_2/c \choose {N-m} } } 
{ {-(\nu_1 + \nu_2)/c  \choose  N} } , 
\end{equation}
where the binomial coefficient ${r \choose n}$ for a real $r$ and
nonnegative integer $k$ denotes
\begin{equation}
\label{eq:binomial-coefficient}
{r \choose k}=\frac{r(r-1)\cdots (r-k+1)}{k!} .
\end{equation}
Note that (\ref{eq:eesnr}) and (\ref{eq:eehy}) are special cases of
(\ref{eq:eepolya}) with $c=0$ and $c=-1$, respectively.  In
(\ref{eq:eehy}) the range of summation can be taken as $m=0,\ldots,N$,
with the convention (\ref{eq:binomial-coefficient}).


\section{Arbitrary discrete distribution with finite support}
\label{sec:arbitrary-distribution}

So far we have discussed how to derive some classical distributions.  We
now show that given any distribution on $\{0,\ldots,N\}$, we can specify a
neutral strategy of Forecaster in a game with $N$ rounds such that Reality
follows the distribution.

Let $q_m\ge 0, m=0,\ldots,N$, $\sum_{m=0}^N q_m=1$, denote an arbitrary
probability distribution on $\{0,\ldots,N\}$.  By decreasing $N$ if
necessary, we assume $q_N > 0$. 
Define
\[
\bar q_m = \frac{q_m +\cdots +q_N}{q_{m-1}+q_m + \cdots + q_N}, \quad m=1,\ldots,N.
\]
Then  
\begin{equation}
\label{eq:qm}
q_m = \bar q_1 \cdots \bar q_m (1-\bar q_{m+1}), \quad m=1,\ldots,N, \ 
  \bar q_{N+1}=0.
\end{equation}
Let
\[
p_n
 = 
\begin{cases} \bar q_n, & \text{if      ${\cS}_{n-1}=n-1$}, \\
                0, & \text{otherwise}.
\end{cases}
\]
The idea here is to let Reality increase ${\cS}_{n-1}$ by 1 with
probability $\bar q_{n-1}$ if ${\cS}_{n-1}=n-1$ or otherwise let him stop at
the current level for the rest of the rounds.  Note that
$p_n=p_n({\cS}_{n-1})$ is indeed a function of ${\cS}_{n-1}$, because it
is written as
\[
p_n=\bar q_n\times I({\cS}_{n-1}=n-1),
\]
where $I(\cdot)$ is the indicator function.  

\medskip
\noindent
\textsc{Biased-Coin Game With Forecaster For Arbitrary Distribution}\\
\textbf{Protocol:}

\parshape=8
\IndentI   \WidthI
\IndentI   \WidthI
\IndentII  \WidthII
\IndentII  \WidthII
\IndentII  \WidthII
\IndentII  \WidthII
\IndentII  \WidthII
\IndentI   \WidthI
\noindent
${\cal K}_0 =\alpha$, $\cS_0=0$, 
$q_m\ge 0, m=0,\ldots,N$, $\sum_{m=0}^N q_m=1$: given\\
FOR $n=1,\dots,N$\\
  Forecaster announces $p_n=\bar q_n \times I({\cS}_{n-1}=n-1)$.\\
  Skeptic announces $M_n\in{\mathbb R}$.\\
  Reality announces $x_n\in \{0, 1\}$.\\
  ${\cal K}_n := {\cal K}_{n-1} + M_n (x_n-p_n)$.\\
  $\cS_n := \cS_{n-1}+x_n$.\\
END FOR\\
\noindent

\begin{figure}
\begin{center}
\setlength{\unitlength}{1mm}
\begin{picture}(50,50)(0,0)
\thicklines
\put(0,0){\vector(1,0){10}}
\put(0,0){\vector(1,1){10}}
\put(10,10){\vector(1,0){10}}
\put(10,10){\vector(1,1){10}}
\put(10,0){\vector(1,0){10}}
\multiput(22,0)(3,0){3}{\circle*{1}}
\multiput(31,0)(3,0){3}{\circle*{1}}
\put(40,0){\vector(1,0){10}}
\multiput(22,10)(3,0){3}{\circle*{1}}
\multiput(31,10)(3,0){3}{\circle*{1}}
\put(40,10){\vector(1,0){10}}
\multiput(22,20)(3,0){3}{\circle*{1}}
\multiput(31,20)(3,0){3}{\circle*{1}}
\multiput(22,22)(3,3){3}{\circle*{1}}
\multiput(31,31)(3,3){3}{\circle*{1}}
\put(40,40){\vector(1,0){10}}
\put(40,40){\vector(1,1){10}}
\put(0,-4){$0$}
\put(8,-4){$1$}
\put(18,-4){$2$}
\put(52,-1){$0$}
\put(52,9){$1$}
\put(52,39){$N-1$}
\put(52,49){$N$}
\put(5,2){\footnotesize$q_0$}
\put(1,5){\footnotesize$\bar q_1$}
\put(13,11){\footnotesize $1 - \bar q_2$}
\put(11,16){\footnotesize $\bar q_2$}
\put(42,41){\footnotesize $1 - \bar q_N$}
\put(41,46){\footnotesize $\bar q_N$}
\end{picture}
\end{center}
\caption{Tree of the game for arbitrary distribution}
\label{fig:1}
\end{figure}
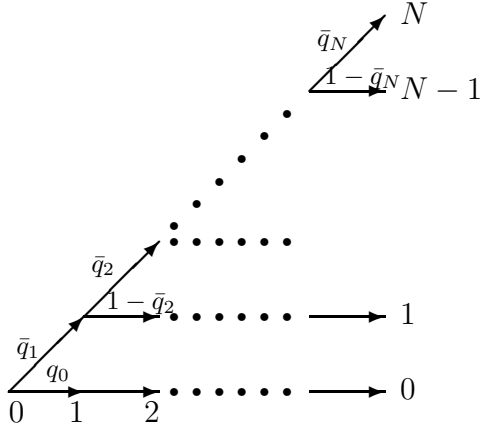

The tree of this game is illustrated in Figure \ref{fig:1}.
For this game we have the following result.

\begin{theorem}
\label{thm:univariate-arbitrary}
  The upper and the lower expected values of $\eta(\cS_N)$ coincide
  and given by
\begin{equation}
\label{eq:eearbitrary}
\bar\EE(\eta(\cS_N)) = \underline\EE(\eta(\cS_N)) = 
\sum_{m=0}^N \eta(m) q_m.
\end{equation}
\end{theorem}

\proof As in the case of hypergeometric distribution 
\[
\bar\eta(0,0)=\sum_{(x_1, \ldots, x_N)\in { \{0,1\}^N}}
\prod_{n=1}^N   p_n(\cS_{n-1})^{x_n} (1-p_n(\cS_{n-1}))^{1-x_n}  
\eta(\cS_N).
\]
In this game, the path leading to ${\cS}_m$ is uniquely determined as
\[
x_1\ldots x_N = 1\ldots 1\, 0 \ldots 0
\]
with $m$ initial 1's.   By (\ref{eq:qm}),  for this path
$\prod_{n=1}^N   p_n(\cS_{n-1})^{x_n} (1-p_n(\cS_{n-1}))^{1-x_n}  = q_m$. 
Therefore $\bar\eta(0,0)=\sum_{m=0}^N \eta(m) q_m$.  The replicating
strategy confirming this candidate price is given as
\[
M_n = 
\begin{cases} 
  \frac{\sum_{m=n}^N \eta(m) q_m}{\sum_{m=n}^N q_m} - \eta({\cS}_{n-1}), 
 & \text{if      ${\cS}_{n-1}=n-1$}, \\
          \infty, & \text{otherwise}.
\end{cases}
\]
\qed

In Theorem \ref{thm:univariate-arbitrary} we have considered a discrete
distribution with the support $\{0,\ldots,N\}$. We can deal with the support 
of the form $\{a,a+1,\ldots,b\}$, by letting $N=b-a$ and setting the
initial value ${\cS}_0=a$.

In Section 8.3 of their book, Shafer and Vovk discuss ``adding tickets''
to make the upper expected value and the lower expected value to
coincide. Note that if Reality's move space has more than two elements in a
single step game, then the upper expected value is generally larger than
the lower expected value.
Theorem \ref{thm:univariate-arbitrary} shows that if we add 
sufficient number of steps to a single step game, the equality of the
upper and the lower prices is achieved.

\section{Multivariate extension}
\label{sec:multivariate}

In the previous sections we have considered univariate random variable
$\cS_N$.  In this section we give a straightforward multivariate extension
of the results of the previous sections.  We employ the multi-label
classification protocol discussed in Vovk, Nouretdinov, Takemura and
Shafer (2005).  


Our extension corresponds to generalizing Binomial distribution to
multinomial distribution.  Let $\cS_N = (\cS^1_N, \ldots, \cS^d_N)$ be a
$d$-dimensional vector. As in multinomial distribution, for the sake of
symmetry, we leave one-dimensional redundancy $N=\cS^1_N+ \cdots +
\cS^d_N$ in the components of $\cS_N$.  Therefore $x_n$ in Rescaled
Biased-Coin Game now corresponds to a 2-dimensional vector $(x_n, 1-x_n)$.
For the general $d$-dimensional case the move space of Reality
\[
{\bf X}=\{e_1, \dots, e_d\} = \{ (1,0,\dots,0), (0,1,0,\dots,0), \dots,
(0,\dots,0,1)\} 
\]
consists of $d$ standard coordinate vectors.

\medskip
\noindent
\textsc{Multilabel Classification Game With Neutral Forecasting Strategy}\\
\textbf{Protocol:}

\parshape=8
\IndentI   \WidthI
\IndentI   \WidthI
\IndentII  \WidthII
\IndentII  \WidthII
\IndentII  \WidthII
\IndentII  \WidthII
\IndentII  \WidthII
\IndentI   \WidthI
\noindent
${\cal K}_0 =\alpha$, $\cS_0=0$ : given\\
FOR $n=1,\dots,N$\\
  Forecaster announces $p_n=p_n(\cS_{n-1}) \in {\mathbb R}^d$.\\
  Skeptic announces $M_n\in{\mathbb R}^d$.\\
  Reality announces $x_n\in {\bf X}$.\\
  ${\cal K}_n := {\cal K}_{n-1} + M_n \cdot (x_n-p_n)$.\\
  $\cS_n := \cS_{n-1}+x_n$.\\
END FOR

\noindent
Here ``$\cdot$'' denotes the standard inner product of ${\mathbb R}^d$.

In the above protocol we took the whole ${\mathbb R}^d$ as the move space
of Forecaster.  Let  
\[
\Delta({\bf X})= \{ (p^1,\ldots,p^d) \mid p^i \ge 0, \sum_{i=1}^d p^i=1 \}
\]
denote the probability simplex spanned by the standard coordinate vectors.
If Forecaster announces $p_n \not\in \Delta({\bf X})$, then by the
hyperplane separation theorem Skeptic can choose $M_n\in {\mathbb R}^d$
such that he becomes infinitely rich immediately, no matter what move
Reality chooses.  
See Vovk, Nouretdinov, Takemura and Shafer (2005) for a discussion of this point.
Therefore we can restrict Forecaster's move space to 
the probability simplex $\Delta({\bf X})$.  Also if $p^i_n=0$ for some $i$,
Skeptic can choose $M^i_n$ arbitrarily large and Reality is forced to
choose $x^i_n=0$.

We also note that there is a redundancy in the move space of Skeptic, once
$p_n$ is restricted to lie in $\Delta({\bf X})$. $M_n +  c (1,\ldots,1)$ for
any $c\in {\mathbb R}$ leads to the same increment of the 
capital process ${\cal K}_n$. However it is often convenient to ignore
this redundancy in specifying Reality's move $M_n$.

For notational simplicity write
\[
p_n^{x_n} = (p^1_n)^{x_n^1}\cdots (p^d_n)^{x_n^d} = p^i_n  \qquad
\mbox{for}\ x_n=e_i.
\]
As a straightforward generalization
of results in the previous sections we have the following theorem.

\begin{theorem}
\label{thm:multivariate}
  The upper and the lower expected values of $\eta(\cS_N)$ coincide
  and given by
\begin{equation}
\label{eq:eearbitrary-multi}
\bar\EE(\eta(\cS_N)) = \underline\EE(\eta(\cS_N)) = 
\sum_{(x_1, \ldots, x_N)\in  {\bf X}^N}
\prod_{n=1}^N   p_n(\cS_{n-1})^{x_n}
\eta(\cS_N).
\end{equation}
\end{theorem}

The line of the proof is the same as in the previous theorems and we omit
the details.   The price $\bar\eta(n,{\cS}_n)$ at time $n$  is defined
recursively by 
\[
\bar\eta(n,{\cS}_n)
= \sum_{i=1}^d p_{n+1}^i({\cS}_n) \bar\eta(n+1,\cS_n + e_i)
\]
and the replicating strategy $M_n^i$ is simply given by
\[
M_n^i =\bar\eta(n,\cS_{n-1} + e_i).
\]
{}From Theorem \ref{thm:multivariate} we can easily derive multinomial
distribution, multivariate hypergeometric distribution as well as
multivariate Polya's distribution.

A generalization of Theorem \ref{thm:univariate-arbitrary} to an arbitrary
$(d-1)$-dimensional discrete distribution of $(\cS_N^1,\ldots.
\cS_N^{d-1})$ with finite support can be explained as follows.  We first
use Theorem \ref{thm:univariate-arbitrary} on the first component
$\cS_n^1$ to derive the one-dimensional marginal distribution of
$\cS_N^1$.  Then, given the realization of the first component, we derive
the conditional distribution of $S_N^2$ given $S_N^1$ by another
application of Theorem \ref{thm:univariate-arbitrary} to the second
component.  We can continue this process up to the $(d-1)${\it th}
component.  The last component $\cS_n^d$ is used as a slack variable.

Finally as an illustration of 
Theorem \ref{thm:multivariate}
we show how the Cox-Ross-Rubinstein formula of Section
\ref{sec:cox-ross-rubinstein} is reduced to our multivariate framework.
Define 
\[
x_n \equiv \frac{\cS_{n} -r\cS_{n-1}}{r^n} .
\]
Furthermore by discounting define
\[
K_n^* \equiv \frac{K_n}{r^n}.
\]
Then the recurrence relation of the Cox-Ross-Rubinstein game
\begin{eqnarray*}
{\cal K}_n 
&=&  {\cal K}_{n-1} + M_n (\cS_n - \cS_{n-1}) +
  (r-1) ({\cal K}_{n-1}- M_n \cS_{n-1})  \\
&=& r {\cal K}_{n-1} + (\cS_n
 - r \cS_{n-1} ) M_n
\end{eqnarray*}
is written as
\[
{\cal K}_n^* = {\cal K}_{n-1}^* + M_n x_n.
\]
Here $x_n$ can take two values
\[
(x^1_n,  x^2_n) 
= \left( \frac{\cS_{n-1} (u-r)}{r^n}, \frac{\cS_{n-1} (d-r)}{r^n} \right).
\]
Rescaling the values we define $d=2$, $x_n^* \in \{e_1, e_2\} = \{(1,0),
(0,1)\}$ and 
\begin{align*}
M_n^* &= (x_n^1, x_n^2), \\
p_n^* &= \left(\frac{-x_n^2}{x_n^1 - x_n^2}, \frac{x_n^1}{x_n^1 - x_n^2}
\right). 
\end{align*}
Then $M_n x_n$ is written as
\[
M_n x_n = M_n^* \cdot (x_n^* - p^*_n).
\]
Therefore the iteration part of the 
the Cox-Ross-Rubinstein game is written as

\medskip
\parshape=5
\IndentI   \WidthI
\IndentII  \WidthII
\IndentII  \WidthII
\IndentII  \WidthII
\IndentI   \WidthI
\noindent
FOR $n=1,\dots,N$\\
  Skeptic announces $M_n^*\in{\mathbb R}^2$.\\
  Reality announces $x_n^*\in \{e_1, e_2\}$.\\
  ${\cal K}_n^*  := {\cal K}^*_{n-1} + M_n^* \cdot (x_n^* - p_n^*)$.\\
END FOR\\

This shows that  the Cox-Ross-Rubinstein formula is also a special case of
our multivariate extension.





\appendix

\section{Preliminaries on game theoretic probability}
\label{sec:preliminries}

Here we summarize preliminary material (Chapter 1 of S\&V)
of the game
theoretic probability.  We also state a basic proposition on the existence
of a replicating strategy and the existence of the game theoretic
expectation in a coherent game.

In this paper all the games are finite-horizon games with $N$ rounds.
Therefore a path of the game is a finite sequence $\xi = x_1 \dots
x_N$ of Reality's moves.  A random variable $x(\xi)$ denotes a payoff to
Skeptic, when Reality chooses the path $\xi$.  Given a strategy ${\cal P}$
of Skeptic, ${\cal K}^{\cal P}$ denotes the capital process of $\cal P$
with zero initial capital.  Furthermore in this paper we only consider
symmetric games, in the sense that if ${\cal P}$ is a strategy of Skeptic,
then $-{\cal P}$ is also a strategy of Skeptic and
\[
{\cal K}^{-{\cal P}}=- {\cal K}^{\cal P} .
\]

The upper expected value $\bar\EE x$  and the lower expected value
$\underline{\EE}x$  of $x$ is defined as
\begin{eqnarray*}
\bar\EE x &=& \inf\{\alpha \mid \exists {\cal P}, \forall \xi, \ 
   {\cal K}_N^{\cal P}(\xi) \ge x(\xi) - \alpha \}, \\
\underline\EE x &=& \sup\{\alpha \mid \exists {\cal P}, \forall \xi, \ 
   {\cal K}_N^{\cal P}(\xi) \ge \alpha - x(\xi) \}.
 \end{eqnarray*}
In a symmetric game $\underline\EE x$ can also  be written as
\begin{equation}
\label{eq:lower1}
\underline\EE x = \sup\{\alpha \mid \exists {\cal P}, \forall \xi, \ 
   \alpha + {\cal K}_N^{\cal P}(\xi) \le x(\xi) \}.
\end{equation}
A game is coherent if Skeptic is not allowed to make money for certain,
i.e., 
\[
\forall {\cal P},\  \exists \xi, \ \KP_N(\xi) < 0.
\]
If a game is coherent, then $\bar\EE x \ge \underline\EE x$ for every
random variable $x$ (Proposition 7.2 of S\&V).  

We call $\cal P$ a {\it replicating strategy} for $x$  with the
replicating initial capital $\alpha \in {\mathbb R}$ if
\[
\alpha + \KP_N(\xi) = x(\xi), \ \forall \xi.
\]
We now state the following basic fact.

\begin{proposition}
\label{prop:replicatioting-strategy}
In a coherent symmetric 
game, suppose that ${\cal P}^*$ is a replicating strategy for
$x$ with the replicating  initial capital $\alpha^*$.  Then
\[
\bar\EE x = \underline\EE x = \alpha^* .
\] 
\end{proposition}

\proof  By definition of $\bar\EE x$ we have $\bar\EE x \le \alpha^*$.
Furthermore in a symmetric game  $\underline\EE x \ge \alpha^*$ follows
from (\ref{eq:lower1}).  Therefore  $\bar\EE x \le \alpha^* \le
\underline\EE x$.  Combining this with the inequality $\bar\EE x \ge
\underline\EE x$ we obtain the proposition. \qed

\bigskip 
We should note that the proof of the inequality $\bar\EE x \ge
\underline\EE x$ in S\&V and the above proof are standard arbitrage
arguments.

\bigskip
\noindent
{\bf Acknowledgment.}\quad  We are grateful for insightful comments by  Vladimir
Vovk, Glenn Shafer, Kei Takeuchi and Masayuki Kumon on earlier drafts of
this paper.

\end{document}